\documentclass[a4paper,12pt]{article}
\begin{document}

\title{Spaces of quasi-invariant measures and convergence in them.}
\author{S.V. Ludkowski.}

\date{05 September 2015}
\maketitle

\begin{abstract}
Spaces of quasi-invariant measures supplied with different
topologies are studied. Their embeddings, projective decompositions,
conditions for their metrizability are investigated. Theorems about
convergence of nets of quasi-invariant measures and their extensions
are proved as well. Moreover, associated with them uniform spaces
are studied.

\footnote{key words and phrases: measure, quasi-invariant, topological space, metric \\
Mathematics Subject Classification 2010: 28A33; 28C15 \\ address: Department of Applied Mathematics, \\
Moscow State Technical University MIREA, \\ av. Vernadsky 78, Moscow
119454, Russia\\
ludkowski@mirea.ru}

\end{abstract}

\section{Introduction.}
Measures and their spaces play very important role in functional
analysis and probability theory \cite{alexmsb941}-\cite{dalshn},
\cite{ferpelruizb}-\cite{vardmsb61}. Among them quasi-invariant
measures are widely used in many branches of mathematics including
harmonic analysis and group representations and geometry. Although
spaces of measures were intensively investigated, but properties of
spaces of quasi-invariant measures remained less studied.
\par Concrete examples of families of quasi-invariant measures are
contained in the cited above literature. They include
quasi-invariant measures on Hilbert spaces and Banach spaces or more
general topological vector spaces relative to proper additive
subgroups
\cite{beldal,dal,kuob}-\cite{ludqupdmb,skorohodb,smolfom,vahtarchobb}.
Then measures quasi-invariant relative to transformation groups
consisting of linear $A$ or non-linear operators $B$ on the
separable real Hilbert space $l_2$ such that $A-I$ or $B'(x)-I$ are
nuclear or more generally Hilbert-Schmidt operators with some other
conditions imposed were described as well, where $I$ denotes the
unit operator, $B'(x)$ notates the strong (Frechet) derivative of
$B$ at $x\in X$. Particularly, they may be Gaussian measures. \par
Also quasi-invariant measures on topological groups which may be
non-locally compact relative to proper subgroups were investigated
in \cite{beldal,dalshn,lujms147:3:08,lusmldg05,lurimut99}-\cite{ludspb}.
More generally, quasi-invariant measures on manifolds and
topological spaces which may be non-locally compact relative to
transformation groups were studied in
\cite{beldal,fidal,lurimut98,ludspb}. For locally compact groups and
topological spaces studies of quasi-invariant measures were begun
earlier (see, for example, \cite{bourm70,hew} and references
therein). Using polyhedral expansions of complete uniform spaces it
is possible to provide abundant families of quasi-invariant measures
on complete uniform spaces relative to their certain transformation
groups \cite{ludqupdmb,ludspb}.\par Apart from spaces of measures,
families of quasi-invariant measures appear to be generally
non-linear, but can be supplied with topological or uniform
structures. This work continues previous publications
\cite{ludqupdmb,lujms147:3:08} of the author on this subject and
treats new aspects of the theory.
\par In those publications convergence of nets of quasi-invariant
measures was studied, but with rather strong conditions on
quasi-invariance factors of measures like a uniform
convergence on each compact subset. In the present work more general
classes of measures are considered and weaker conditions on
quasi-invariance factors and topologies are imposed. Furthermore,
new general properties of topological spaces of quasi-invariant
measures are investigated besides convergence of measure nets.
\par In this article spaces of quasi-invariant measures supplied with
different topologies are studied. Their embeddings, projective
decompositions, conditions for their metrizability are investigated.
Theorems about convergence of nets of quasi-invariant measures and
their extensions are proved as well. Moreover, associated with them
uniform spaces are studied.
\par As it is known under definite conditions the Radon-Nikodym derivative
of one $\sigma $-smooth measure relative to another may exist \cite{bogachmtb,bourm70}.
Nevertheless in the present work it is not supposed in advance that a quasi-invariant
measure should have a quasi-invariance factor, that is the Radon-Nikodym derivative
of a transformed measure relative to its initial measure.
\par In Theorem 2 two different topologies on spaces of measures
are compared taking into account transformation groups. In
Proposition 5 a relation between quasi-invariant measures with the
Souslin number of a topological space is elucidated. Theorem 6
reveals an interplay between bounded quasi-invariant functionals and
$\sigma $-smooth quasi-invariant measures. An equivalence relation
on a space of continuous functions induced by a quasi-invariance
factor is studied (see Proposition 7 and Corollary 8). \par It is
shown in Proposition 9 that in a topological space of
quasi-invariant measures there are closed linear subspaces relative
to a weak topology. Theorem 11 about extensions of quasi-invariant
measures from Tychonoff topological spaces onto their Wallman
extensions (compactifications) is proved. Metrizability of spaces of
quasi-invariant measures under definite conditions is described in
Theorem 12. Some results on ranges of quasi-invariant measures are
proved in Lemma 13 and Proposition 14. Conditions for a space of
quasi-invariant measures to be dense or not in a space of measures
in a weak topology are investigated (see Theorem 15 and Proposition
16). \par Uniform spaces associated with algebras of sets and spaces
of quasi-invariant measures and actions on them of groups are
studied (see Theorems 18 and 19). Then Theorem 21 is proved about
convergence of sequences of quasi-invariant measures under rather
mild conditions. Extensions of quasi-invariant measures for
topological groups are investigated as well (see Theorem 22).
Decompositions of quasi-invariant measures and their spaces with the
help of inverse mapping systems are studied in Theorems 23 and 24.
\par All main results of this paper are obtained for the first time.
They can be used for further studies of spaces of quasi-invariant measures,
groups algebras, representations of groups and algebras.

\section{Families of quasi-invariant measures.}
\par To avoid misunderstandings we first remind some definitions which may vary
in the literature.
\par {\bf 1. Definitions.} Let $X$ be a $T_1\cap T_{3.5}$ topological space
and let $G$ be a group acting on $X$ so that to each element $g\in
G$ there corresponds a mapping $h_g: X\to X$ satisfying the
conditions $h_{g}\circ h_{j} = h_{gj}$ for all $g, j\in G$ and
$h_e=id$, where $id(x)=x$ for each $x\in X$ designs the identity
mapping, $e$ denotes the neutral element in $G$. We shortly denote
$h_g(x)=gx$ for each $g\in G$ and $x\in X$.
\par Let $C(X,{\bf F})$ (or $C_b(X,{\bf F})$) be the space of all
continuous (or continuous and bounded respectively) functions from
the topological space $X$ into the field $\bf F$, where the field
$\bf F$ is either real ${\bf F}={\bf R}$ or complex ${\bf F}={\bf
C}$. Put ${\sf Z} := \{ Z: Z=f^{-1}(0), f \in C_b(X,{\bf R}) \} $
and ${\sf U} := \{ U: U= X\setminus Z, ~ Z \in {\sf Z} \} $, also
the notation ${\cal F}={\cal F}(X)$ will be used for the minimal
algebra containing these families $\sf Z$ and $\sf U$, then ${\cal
B}={\cal B}(X)$ will stand for the minimal $\sigma $-algebra
containing $\sf Z$ and $\sf U$. \par Henceforth, it is supposed that
$g: {\cal F}\to {\cal F}$ and $gX=X$ for each $g\in G$, if something
other will not be outlined.
\par Measures $m, n: {\cal F} \to \bf  F$ are called equivalent $m\sim n$ if
$|m| <<|n|$ and $|n|<<|m|$, where $|m|$
 denotes the variation of the measure $m$, while the notation $|m| <<|n|$
means that $|m|$ is absolutely continuous relative to $|n|$. \par A
measure $m: {\cal F}\to \bf F$ is called (left) quasi-invariant
relative to a group $G$ if $m^g$ is equivalent to $m$ for each $g\in
G$, where $m^g(A) := m(g^{-1}A)$ for each $A\in \cal F$. Then
$d_m(g,x) := m^g(dx)/m(dx)$ denotes the (left) quasi-invariance
factor of $m$ wherever it exists, where $g\in G$, $x\in X$ and it is
supposed that $d_m(g,x)\in {\bf F}$.
\par  We shall use the notation $M(X,{\bf F}, {\cal F})$ for the family
of all measures on $\cal F$ with values in $\bf F$. The
corresponding family of all quasi-invariant relative to the group
$G$ measures on the measurable space $(X,{\cal F})$ is denoted by
$Q(X,{\bf F},{\cal F},G)$ and its subfamily of non-negative measures
by $Q^+(X,{\bf R},{\cal F},G)$.
\par Let a function $d: G\times X\to {\bf F}$
satisfy the conditions: \par $(1)$ $d(e,x)=1$ for each $x\in X$,
where $e\in G$ denotes the neutral element of the group $G$,
\par $(2)$ $d(g,x)\ne 0$ for every $g\in G$ and $x\in X$, also let
\par $(3)$ $d$ satisfy the co-cycle condition:
$d(sg,x) = d(g,s^{-1}x)d(s,x)$ for all $s, g \in G$ and $x\in X$.
\par We denote by $Q^d(X,{\bf F},{\cal F},G)$
the family of all measures $m \in M(X,{\bf F},{\cal F})$  such that
the measure $m$ is left quasi-invariant with existing
quasi-invariance factor $d_m(g,x) = d(g,x)$ for each $g\in G$ for
$m$-almost all $x\in X$. Its subfamily of non-negative measures is
denoted by $Q^{+,d}(X,{\bf R},{\cal F},G)$. If some data like $\bf
F$ or $G$ or $\cal F$ are specified, they will be omitted for
shortening of the notation.

\par {\bf 2. Theorem.} {\it The families of neighborhoods
$$(1)\quad N_s(m_0;f_1,...,f_n;G;y) := $$ $$\{ m: m\in M(X), ~ \forall
r=1,...,n ~ sup_{g\in G} |\int_X f_r(x)(m^{g}(dx) - m^{g}_0(dx)|<y
\} ;$$
$$(2)\quad N_w(m_0;f_1,...,f_n;g(1),...,g(k);y) := $$ $$\{ m: m\in M(X),
~ \forall r=1,...,n, j=1,...,k ~ |\int_X f_r(x)(m^{g(j)}(dx) -
m^{g(j)}_0(dx)|<y \}
$$
induce $T_1\cap T_{3.5}$ topologies $\tau _s(G)$ and $\tau _w(G)$ on
$M(X)=M(X,{\bf F},{\cal F})$ and $Q(X)\hookrightarrow M(X)$, where
$f_1\in C_b(X,{\bf F}),...,f_n\in C_b(X,{\bf F})$, $~ n, k\in {\bf
N}$, $~y>0$ and $g(1)\in G,...,g(k)\in G$, $~Q(X)=Q(X,{\bf F},{\cal
F},G)$. The first topology is generally stronger than the second
one, when $X$ and $G$ are infinite. If $h_g: X\to X$ is the
homeomorphism for each $g\in G$, then the topology $\tau _w(G)$ is
equivalent with the usual weak topology $\tau _w = \tau _w( \{ e \}
)$.}
\par {\bf Proof.} If $m\in M(X)$, then $m^g\in M(X)$ as well, since
$g: {\cal F}\to {\cal F}$, where $m^g(E):=m(g^{-1}E)$ for each $E\in
\cal F$ and $g\in G$, though $m^g$ need not be equivalent with $m$.
At first we verify from Formulas $(1)$ and $(2)$, that
\par $N_s(m_0;f_1,...,f_n;G;y)\cap N_s(m_0;h_1,...,h_k;G;y)=
N_s(m_0;f_1,...,f_n,h_1,...,h_k;G;y)$ and
\par $N_w(m_0;f_1,...,f_n;g(1),...,g(k);y)\cap
N_w(m_0;h_1,...,h_l;u(1),...,u(t);y)=$\par $
N_w(m_0;f_1,...,f_n,h_1,...,h_l;g(1),...,g(k),u(1),...,u(t);y)$ \\
for every functions $f_1\in C_b(X,{\bf F}),..., f_n\in C_b(X,{\bf
F}), h_1\in C_b(X,{\bf F}),..., h_l\in C_b(X,{\bf F})$ and group
elements $g(1)\in G,..., g(k)\in G, u(1)\in G,..., u(t)\in G$. Also
for each $m_0\in M(X)$ there exist non-void neighborhoods
$N_s(m_0;f_1,...,f_n;G;y)$ and
$N_w(m_0;f_1,...,f_n;g(1),...,g(k);y)$. Therefore, the families
${\sf B}_s(G)$ and ${\sf B}_w(G)$ given by formulas $(1)$ and $(2)$
respectively are bases of topologies, since they satisfy Conditions
1.1$(B1,B2)$ \cite{eng}, hence they induce topologies, which were
denoted by $\tau _s(G)$ and $\tau _w(G)$ correspondingly above.
\par In view of Theorem II.1 \cite{vardmsb61} the topological space
$(M(X),\tau _w(G))$ is completely regular (Tychonoff), that is
$(M(X),\tau _w(G))\in T_1\cap T_{3.5}$. Since $Q(X)\subset M(X)$,
then the topology $\tau _w(G)$ on $M(X)$ induces the corresponding
topology on $Q(X)$ and hence $(Q(X),\tau _w(G))\in T_1\cap T_{3.5}$.
\par If the topological space $X$ and the group $G$ are infinite,
then the topology $\tau _s(G)$ is generally stronger, than $\tau
_w$, since the neighborhood $N_s(m_0;f_1,...,f_n;G;y)$ with
nonconstant continuous bound3ed functions $f_1,...,f_n$ can not be
obtained in general as the finite intersection of weak neighborhoods
\par $N_w(m_0;f_1,...,f_n;g(1),...,g(k);y)$, where $n, k\in {\bf N}$.
\par If $m_0\ne m_1 \in M(X)$, then there exists $f\in C_b(X,{\bf F})$
so that $m_0(f)\ne m_1(f)$ (see also \cite{alexmsb941,bogachmtb}),
where
$$m(f) := \int_X f(x)m(dx),$$ since $X\in T_1\cap T_{3.5}$. Then
$$\sup_{g\in G} |m_0^g(f)-m_1^g(f)|\ge |\int_X
f(x)[m_0(dx)-m_1(dx)]|=y >0,$$ since $m^e(dx)=m(dx)$ for the neutral
element $e$ of the group $G$. Therefore, $N_s(m_0;f;G;y/4)\cap
N_s(m_1;f;G;y/4) = \emptyset $, consequently, $(M(X),\tau _s(G))\in
T_2$ and hence $(M(X),\tau _s(G))\in T_1$ (see also \S 1.5
\cite{eng}).
\par Let now $m_0\in M(X)$ and let $J$ be a closed subset in the topological space
$(M(X), \tau _s(G))$. We take the neighborhood \par $U :=
N(m_0;f_1,...,f_n;G;y)$ \\ of a measure $m_0\in M(X)$ and mention
that $$U =\bigcap_{j=1}^n N(m_0;f_j;G;y),$$ where $y>0$. Put
$$u_j(m) := \min (\sup_{g\in G} |m^g(f_j)-m_0^g(f_j)|/y, 1) ,$$ then
$u_j(m)$ is continuous relative to the $\tau _s(G)$ topology on
$M(X)$, where $y>0$ is a marked number. Moreover, $u_j(m_0)=0$, also
$u_j(m)=1$ on $M(X) - N(m_0;f_1,...,f_n;G;y)$. Thus the function
$u(m) := \max (u_1,....,u_n)$ is such that $u(m_0)=0$ and $u(m)=1$
on $M(X)-U$. This means that $(M(X),\tau _s(G))\in T_{3.5}$. \par
Since $Q(X)\subset M(X)$, the topology $\tau _s(G)$ on $M(X)$
induces it on $Q(X)$ and consequently, $(Q(X),\tau _s(G))\in T_1\cap
T_{3.5}$ as well.
\par Let now $h_g: X\to X$ be the homeomorphism for each $g\in G$, then
to each bounded continuous function $f\in C_b(X,{\bf F})$ there
corresponds $f^g\in C_b(X,{\bf F})$, where $f^g(x) := f(g^{-1}x)$
for each $x\in X$. Therefore, from the equality
$$\int_X f(x)m^{g}(dx) = \int_X f^{g^{-1}}(y)m(dy)$$ for each $f\in
C_b(X,{\bf F})$ and $g\in G$ it follows that $$(3)\quad
N_w(m_0;f_1,...,f_n;g(1),...,g(k);y) := $$
$$\bigcap_{i=1}^n\bigcap_{j=1}^k N_w(m_0;f_{i,j};y)
,$$ where $f^{g}(y):= f(g^{-1}y)$, $~f_{i,j}=f_i^{g(j)^{-1}}$ for
each $i=1,...,n$ and $j=1,...,k$. Thus the topology $\tau _w(G)$ on
$M(X)$ is equivalent to the weak topology $\tau _w = \tau _w( \{ e
\} )$, since their bases ${\sf B}_w(G)$ and ${\sf B}_w( \{ e \} )$
are equivalent in the considered case.

\par {\bf 3. Corollary.} {\it $(Q^+(X),\tau _w)$ is closed in the
topological space $(Q(X),\tau _w)$ and $(Q^+(X),\tau _s(G))$ is
closed in the topological space $(Q(X),\tau _s(G))$.}
\par {\bf Proof.} If $(m_j: j\in J)$ is a net of non-negative
measures converging to $m$ in $(Q(X),\tau _w)$ or in $(Q(X),\tau
_s(G))$ respectively, then its limit $m$ is also non-negative.

\par {\bf 4. Definitions.} A subset $V$ in the space $C(X,{\bf R})$ is called
bounded if two functions $f\in C(X,{\bf R})$ and $h \in C(X,{\bf
R})$ exist such that $f\le u \le h$ for each $u\in V$.
\par It is said that a linear functional $p: C(X,{\bf R}) \to \bf R$ is
bounded, when $p(V)$ is bounded in $\bf R$ for each subset $V$
bounded in $C(X,{\bf R})$.
\par A bounded linear functional $p$ is called $\sigma $-smooth
(or $\tau $-smooth) if \par $\lim_np(f_n)=0$ for each  sequence
(net) $(f_n: n\in A)$ in $C(X,{\bf R})$ so that $f_n\downarrow 0$ \\
(that is, $f_n(x)\ge f_m(x)\ge 0$ for each $m\ge n\in A$ and $lim_n
f_n(x)=0$ for each $x\in X$), where $A=\bf N$ (or $A$ is a directed
set respectively). Then a $\bf C$-linear functional $p: C(X,{\bf
C})\to C(X,{\bf C})$ is bounded, if $p=p_1+ip_2$ and their
restrictions $p_j|_{C(X,{\bf R})} : C(X,{\bf R})\to \bf R$ are $\bf
R$-linear and bounded for $j=1$ and $j=2$, where $p_j : C(X,{\bf
C})\to C(X,{\bf C})$ are $\bf C$-linear functionals for $j=1$ and
$j=2$; $~i=\sqrt{-1}$.
\par The corresponding measures on $\cal B$ are called $\sigma $-smooth or $\tau
$-smooth respectively, their families are denoted by $M_{\sigma
}(X,{\bf F},{\cal F})$ or $M_{\tau }(X,{\bf F},{\cal F})$, while the
corresponding families of quasi-invariant relative to $G$ measures
will be denoted by $Q_{\sigma }(X,{\bf F},{\cal F},G)$ or $Q_{\tau
}(X,{\bf F},{\cal F},G)$. When some data like $\bf F$, $G$ or $\cal
F$ are outlined, they may be omitted to make notation shorter.

\par {\bf 5. Proposition.} {\it If a Hausdorff topological group $Y$ has a left
quasi-invariant $\sigma $-finite $\sigma $-smooth measure $m: {\cal
B}(Y)\to [0,\infty ]$ relative to a dense subgroup $G$, also
$0<m(U)$ for each open symmetric subset $U=U^{-1}$ in $Y$, then the
Souslin number $s(Y)$ of $Y$ is countable, $s(Y)\le \aleph _0$.}
\par {\bf Proof.} The measure $m$ is quasi-invariant, hence if $m(V)=0$
for some subset $V\in {\cal B}(Y)$ in $Y$, then $m^g(V)=0$ for all
$g\in G$.
\par On the other hand, for each symmetric neighborhood $W=W^{-1}$ of
the neutral element $e$ in $Y$ and each open subset $T$ in $Y$ an
element $g\in G$ exists so that $g^{-1}W\cap T\ne \emptyset $. The
measure $m$ is $\sigma $-finite, that is a countable disjoint family
of subsets $Y_j\in {\cal B}(Y)$ exists such that $0<m(Y_j)<\infty $
and $Y=\bigcup_{j=1}^{\infty }Y_j$, where $Y_j\cap Y_k=\emptyset $
for each $j\ne k$. Therefore, for each open symmetric subset
$W=W^{-1}$ in $Y$ we have that a number $j$ exists for which
$m(W\cap Y_j)>0$, since $0<m(W)$. Hence $m^g(W\cap Y_j)>0$ for each
$g\in G$ and consequently, $m(g^{-1}W)>0$. \par Suppose the contrary
that $s(Y)>\aleph _0$, then a family $W_b=W_b^{-1}$ of open subsets
and $g_b\in G$ of elements would exist so that $g_b^{-1}W_b\cap
g_c^{-1}W_c=\emptyset $ for each $b\ne c \in J$, where $J$ is a set
of the cardinality $card (J)>\aleph _0$. Then $card \{ (b,i): ~ b\in
J, i\in {\bf N}, ~ m(g_b^{-1}W_b\cap Y_i)>0 \}
>\aleph _0$. This implies the existence of a positive rational number $y$ so
that $card \{ (b,i): ~ b\in J, i\in {\bf N}, ~ m(g_b^{-1}W_b\cap
Y_i)> y \} >\aleph _0$ and consequently, there exists $j$ so that
$card \{ b: ~ b\in J, ~ m(g_b^{-1}W_b\cap Y_j)> y \}
>\aleph _0$, since $\aleph _0\aleph _0=
\aleph _0$. But this leads to the contradiction, since the measure
$m$ is $\sigma $-smooth and $\sigma $-finite, $0<m(Y_j)<\infty $.
\par We remind the following definitions.

\par {\bf 6. Theorem.} {\it Let $p$ be a bounded linear functional on
$C(X,{\bf F})$ such that
\par $(1)$ $p(f^{g^{-1}}(x))=p(d(g,x)f(x))$ \\ for each $f\in C_b(X,{\bf F})$ and $g\in G$, \\ where
$d(g,x)$ is a continuous function in the $x$ variable for each $g\in
G$, where $x\in X$, $f^g(x):=f(g^{-1}x)$. Then a quasi-invariant
$\sigma $-smooth measure $m$ exists $m\in Q^d_{\sigma }(X,{\bf F})$
so that
$$(2)\quad p(f) = \int_X f(x) m(dx)$$ for each $f\in C(X,{\bf R})$.
Moreover, $d$ satisfies Conditions 1$(1-3)$ (see Definitions 1)
$m$-almost everywhere on $X$.}
\par {\bf Proof.} In view of Theorem I.23 \cite{vardmsb61} and \S 4
there exists $m\in M_{\sigma }(X,{\bf F})$ such that Formula $(2)$
is fulfilled for each $f\in C(X,{\bf F})$. Mention that each ball
\par $B(f,r) := \{ u: u \in C_b(X,{\bf R}); \| u \| \le r \} $ \\ of
radius $0<r<\infty $ in $C_b(X,{\bf R})$ is bounded in $C(X,{\bf
R})$, since $ -r -s\le u(x)\le r+s$ for each $x \in X$, where $$s=
\| f \| := \sup_{x\in X} |f(x)|,$$  while $1_X\in C_b(X,{\bf
R})\subset C(X,{\bf R})$, where $1_X(x)=1$ for each $x\in X$. By the
conditions of this theorem the functional $p$ is bounded and linear
on $C(X,{\bf F})$, consequently, its restriction on $C_b(X,{\bf F})$
is continuous as well according to Definition 4, since $C_b(X,{\bf
C})=C_b(X,{\bf R})\oplus_{\bf R}iC_b(X,{\bf R})$. On the other hand,
we have that
$$(3)\quad \int_X f(x)m^{g}(dx) = \int_X
f^{g^{-1}}(y)m(dy)$$ for each $f\in C_b(X,{\bf F})$ and $g\in G$,
since $g: {\cal F}\to {\cal F}$ and $gX=X$, where
$m^{g}(dx):=m(g^{-1}dx)$. From Formulas $(1)$ and $(3)$ it follows
that
$$p(f^{g^{-1}}) = \int_X f(x)m^{g}(dx)$$ and $m^{g}\in  M_{\sigma }(X,{\bf
F})$ for all $g\in G$ and $f\in C_b(X,{\bf F})$. Moreover, Theorems
I.5, I.22 \cite{vardmsb61} and Equalities $(1,3)$ imply that $m^g$
is equivalent with $m$ for each $g\in G$, since $f\in C_b(X,{\bf
F})$ is arbitrary, consequently, $m\in Q^d_{\sigma }(X,{\bf R})$.
Then from Theorems I.5, I.22 \cite{vardmsb61} and $(1,2)$ we deduce
that the identity $m^g(dx)/m(dx)= d(g,x)$ is valid $m$ almost
everywhere on $X$ for each $g\in G$. Indeed, as the function of the
$x$ variable $d(g,x)$ belongs to $C(X,{\bf F})$ and hence
$d(g,)C_b(X,{\bf F})\subset C(X,{\bf F})$ for each $g\in G$, where
$m^g(dx)/m(dx)$ denotes the Radon-Nikodym derivative (see, for
example, \cite{bogachmtb,bourm70}). Then
$$m^{sg}(dx)/m(dx)=[m(g^{-1}s^{-1}dx)/m(s^{-1}dx)] [
m(s^{-1}dx)/m(dx)]$$ $$ \Rightarrow d(sg,x) = d(g,s^{-1}x)d(s,x)$$
for each $s, g \in G$ and $m$-almost everywhere in the $x$-variable,
$x\in X$. Furthermore, $m^e(dx)/m(dx)=1$ for each $x\in X$, since
$m^e=m$. Therefore, the function $d$ satisfies Conditions 1$(1-3)$
(see Definitions 1) $m$-almost everywhere on $X$ in the $x$ variable
and for each marked elements $s, g \in G$ in Conditions 1$(2,3)$.

\par {\bf 7. Proposition.} {\it Suppose that a function $d: G\times X\to \bf F$
satisfies Conditions 1$(1-3)$ (see Definitions 1) and is continuous
in the $x$ variable on $X$ for each $g\in G$ and a mapping $g: X\to
X$ is continuous for each element $g\in G$.
\par Then there exists an equivalence relation on
$C(X, {\bf F})$ induced by this function $d$.}
\par {\bf Proof.} For any functions $f\in C(X,{\bf F}), h \in C(X,{\bf F})$ we shall say
that they are equivalent $f \Upsilon _d h$ if and only if there
exists an element $g\in G$ so that $f(x) = d(g,x)h(g^{-1}x)$ for all
$x\in X$. Therefore, if $h\in C(X,{\bf F})$, then
$h(g^{-1}x)=h^g(x)\in C(X,{\bf F})$ as the functions of the $x$
variable for each marked element $g\in G$ and hence
$d(g,x)h(g^{-1}x) \in C(X,{\bf F})$, since $d(g,x)$ is continuous in
the $x$ variable and the mapping $g^{-1}: X\to X$ is continuous.
\par Then in virtue of Condition 1$(1)$ (see Definitions 1) we get that $f \Upsilon _d f$, since $f(x)
= d(e,x)f(e^{-1}x)$ for all $x\in X$. That is the relation $
\Upsilon _d $ is reflexive.
\par If $f\Upsilon _d h$ and $h\Upsilon _d u$, then there exist $t,
s\in G$ so that $f(x) = d(t,x)h(t^{-1}x)$ and $h(x) =
d(s,x)u(s^{-1}x)$ for each $x\in X$, consequently, \par $f(x) =
d(t,x)d(s,t^{-1}x) u(s^{-1}t^{-1}x) = d(ts,x) u((ts)^{-1}x)$ \\ due
to the cocycle Condition 1$(3)$. Thus $f \Upsilon _d u$ and the
relation $\Upsilon _d$ is transitive.
\par If $f\Upsilon _d h$, then $h(y) = h(t^{-1}x) = f(x)/ d(t,x)= d(s,y) f(s^{-1}y) $
according to Conditions 1$(1,3)$, where $x=ty$ and $s=t^{-1}$.
Therefore, the relation $\Upsilon _d$ is symmetric. Thus $\Upsilon
_d$ is the equivalence relation on $C(X, {\bf F})$.
\par {\bf 8. Corollary.} {\it Let $m$ be a measure $m\in M(X,{\bf F})$,
let a mapping $g: X\to X$ be continuous and let $d(g,x)$ be a
function continuous in the $x\in X$ variable for each $g\in G$ and
satisfying Conditions 1$(1-3)$ (see Definitions 1). Then the
following conditions are equivalent:
\par $(1)$ for each $f\in C_b(X,{\bf F})$ and $h\in C(X,{\bf F})$ so that $f\Upsilon
_dh$ the equality is fulfilled $m(f)=m(h)$, where
$$m(f) := \int_X f(x)m(dx);$$
\par $(2)$ $m\in Q^d(X,{\bf F})$. }
\par {\bf Proof.} Suppose that Condition $(1)$ of this corollary is fulfilled and $f\Upsilon
_dh$, where $f\in C_b(X,{\bf R})$. In view of Proposition 7 an
element $g\in G$ exists such that $h(x) = d(g,x) f(g^{-1}x) $ for
each $x\in X$, since the mapping $g^{-1}: X\to X$ is continuous. The
change of the variable gives the equality
$$ \int_X u(gy)m(dy) =\int_X u(x)m^g(dx)$$ for each $u\in C_b(X,{\bf F})$,
since $g: {\cal F}\to {\cal F}$ and $gX=X$. The continuous bounded
function $f\in C_b(X,{\bf R})$ is arbitrary, consequently,
$m^g(dx)/m(dx)=d(g,x)$ almost everywhere on $X$ relative to $m$ for
each $g\in G$, since
$$m(f) = \int_X f(y)m(dy)\mbox{ and }m(h) = \int_X
f(g^{-1}x)d(g,x)m(dx).$$ That is $m\in Q^d(X,{\bf F})$.
\par If Condition $(2)$ is satisfied and $f\in C_b(X,{\bf F})$ and $h\in C(X,{\bf F})$
so that $f\Upsilon _dh$ then by Proposition 7 there exists $g\in G$
for which $h(x) = d(g,x) f(g^{-1}x) $ for each $x\in X$, then
$$m(h) = \int_X f(g^{-1}x) d(g,x)m(dx)= \int_X f(g^{-1}x) m(g^{-1}dx)= \int_X f(y)m(dy) ,$$
since $g: {\cal F}\to {\cal F}$ and $gX=X$, consequently, $m(f) =
m(h) $.

\par {\bf 9. Proposition.} {\it Let $d(g,x)$
be a function continuous in the $x\in X$ variable for each $g\in G$
and satisfying Conditions 1$(1-3)$ (see Definitions 1). Then
$Q^d(X,{\bf F})$ is a closed linear subspace in $M(X,{\bf F})$
relative to the $\tau _w$ topology.}
\par {\bf Proof.} If $m_1, m_2\in Q^d(X,{\bf F})$ and $ ~ a, b \in {\bf F}$, then for
the measure $m(dx)=am_1(dx)+bm_2(dx)$ the identities
$m^g(dx)=am_1^g(dx)+bm_2^g(dx)=
ad(g,x)m_1(dx)+bd(g,x)m_2(dx)=d(g,x)m(dx)$ are valid for each $g\in
G$. Therefore, $Q^d(X,{\bf F})$ is the linear space over $\bf F$.
\par Consider an arbitrary net  $(m_k: k\in K)$ in $Q^d(X,{\bf F})$
converging in $M(X,{\bf F})$ to some measure $m$ relative to the
$\tau _w$ topology, where $K$ is a directed set. Then we infer that
$$\int_Xf(x)m(g^{-1}dx) = \lim_k \int_Xf(x)m_k^g(dx) $$
$$=  \lim_k \int_X f(x)d(g,x)m_k(dx) = \int_X f(x)d(g,x)m(dx)$$ for each $f\in C_b(X,{\bf
F})$ and $g\in G$, since $g: {\cal F}\to {\cal F}$ and $gX=X$.
Therefore, there exists the Radon-Nikodym derivative
$m(g^{-1}dx)/m(dx) = m^g(dx)/m(dx) =d(g,x)$ almost everywhere
relative to the measure $m$ for each element $g\in G$. Thus $m\in
Q^d(X,{\bf F})$.

\par {\bf 10. Remark.} On the other hand, the set $Q(X,{\bf F})$ is not a
linear space, when $X$ and $G$ are non-trivial, because different
measures may have different quasi-invariance factors. Mention that
generally $Q(X,{\bf F})$ is not closed in $(M(X,{\bf F}), \tau _w)$,
when $X$ and $G$ are non-trivial. For example, consider a sequence
of Gaussian measures $\lambda _k$ on the real field $\bf R$, that is
on the measurable space $({\cal B}({\bf R}),{\bf R})$, with the same
mean value $v$ and with dispersion $D_k$ tending to zero. Each
measure in this sequence is equivalent to the Lebesgue measure and
is quasi-invariant relative to the additive group $({\bf R},+)$ of
the real field. But this sequence $\lambda _k$ converges in $(M({\bf
R},{\bf R}), \tau _w)$ to the Dirac measure $\delta _v$ with the
support in $v$, which is certainly not quasi-invariant relative to
$({\bf R},+)$.
\par Analogous examples can be considered for the Euclidean space
${\bf R}^n$ relative to its additive group $G=({\bf R}^n,+)$, $n\in
\bf N$, the separable Hilbert space $l_2$ with $G$ being a dense
proper subgroup of the additive group $(l_2,+)$ and for a separable
Banach space $X$ with a dense proper subgroup of the additive group
$(X,+)$ (see about Gaussian measures and their generalizations on
Hilbert spaces and Banach spaces in
\cite{dal,kuob,skorohodb,vahtarchobb,ludanmat2002,ludqupdmb}).
\par By $wX$ is denoted the Wallman extension of a topological space
$X$, while $cl_A B$ denotes
the closure of a subset $B$ in a topological space $A$.
\par Let ${\cal D}(X)$ be the family of all closed subsets in $X$
and let ${\cal U}(X)$ (and ${\cal U}_0(X)$) be the family of all
ultrafilters (of all vanishing ultrafilters correspondingly) in
${\cal D}(X)$, where by definition an ultrafilter $Y$ is vanishing,
if $\cap \{ A: ~ A\in Y \} = \emptyset $.
\par A measure $m: {\cal F}\to \bf F$ is called real, if $$\lim_{ \{ F
\} }m(F)=0$$ for each vanishing ultrafilter $ \{ F \} $ in ${\cal
D}(X)$.

\par {\bf 11. Theorem.} {\it Suppose that a group $G$ acts on
a topological space $X$ continuously $h_g: X\to X$ for each $g\in
G$. Then there exists a mapping $v: Q(X,{\bf C})\to Q(wX,{\bf C})$
so that
\par $(1)$ $p(F') = \lim _{ \{ F \} } m(F)$ \\ for each $m \in Q(X,{\bf
R})$ with $p=v(m)$ and each closed subset $F'$ in $wX$, where $ \{ F
\} $ is the ultrafilter in $X$ satisfying the condition $\cap \{
cl_{wX} F \} = F'$; particularly, \par $p(cl_{wX}F)=m(F)$.
\\ Moreover, $m$ is real if and only if $p(F')=0$ for all closed $F'$
in $wX$ satisfying the condition $F'\subset wX\setminus X$. If
additionally $m$ and $v(m)$ are $\sigma $-smooth, also the variation
of $m$ is finite on $X$, then
\par $(2)$ $d_{v(m)}(g,x)$ exists for every $m \in Q(X,{\bf C})$, $g\in G$ and $x\in
wX$, also $d_{v(m)}(g,x)=d_m(g,x)$ for each $g\in G$ and $m$-almost
everywhere on $X$.}
\par {\bf Proof.} The Wallman extension $wX$ of $X$ is provided
by adjoining to $X$ new points which are the vanishing ultrafilters
in $X$. To each closed subset $F$ in $X$ is posed a closed set $\bar
F$ in $wX$ by adjoining to it all vanishing ultrafilters into which
enters $F$. Then the intersection $\bigcap_k {\bar F}_k$ of any
number of such ${\bar F}_k$ having $\bigcap_k {\bar F}_k\cap X$
closed in $X$ is considered as closed in $wX$.
\par Since the mapping $h_g: X\to X$ is continuous and $h_{g^{-1}}\circ h_g=h_e=id$ and
$gX=X$ for each $g\in G$, then $h_g: X\to X$ is the homeomorphism.
Therefore, if $Y = \{ F \} $ is an ultrafilter in $X$, then $gY = \{
gF: ~ F\in Y \} $ is also the ultrafilter in $X$ for each $g\in G$,
since $h_g: X\to X$ is the homeomorphism, where $gA= \{ gx: ~ x\in A
\} $, $ ~ gx=h_g(x)$. Moreover, if $Y$ is maximal, then $gY$ is
maximal as well. On the other hand, $\bigcap \{ F: ~ F\in Y \} =
\emptyset $ if and only if $\bigcap \{ gF: ~ F\in Y \} = \emptyset
$. Therefore, if $\bar F$ is closed in $wX$, then $g\bar F$ is
defined and also closed in $wX$ for each $g\in G$ and hence
$g\bigcap_k{\bar F}_k = \bigcap_k g {\bar F}_k$. This implies that
if $V$ is closed in $wX$, then $gV$ is closed in $wX$ for each $g\in
G$. Thus the homeomorphism $h_g: X\to X$ has the homeomorphic
extension $h_g: wX\to wX$ for each $g\in G$, since $G$ is the group
and $h_{g^{-1}}\circ h_g=h_e=id$ is the identity mapping.
\par Take an arbitrary quasi-invariant measure $m \in Q(X,{\bf C})$
and an ultrafilter $ \{ F \} $ in $X$. Since $Q(X,{\bf C})\subset
M(X,{\bf C})$, there exists $p=v(m)\in M(wX,{\bf C})$ so that
$p(F')=\lim_{ \{ F \} } m(F)$ and $m$ is real if and only if
$p(F')=0$ for all closed $F'$ in $wX$ satisfying the condition
$F'\subset wX\setminus X$ according to Theorem 12.3
\cite{alexmsb941}.
\par It remains to prove that $p$ is quasi-invariant on $wX$. We infer that
$$(3)\quad \lim_{ \{ F \} } m (g^{-1}F) =  \lim_{ \{ F \} } m^g (F)
= p^g(F')=p (g^{-1}F') =\int_{F'} p(g^{-1}dy)$$ for each $F'$ closed
in $wX$. Therefore, from $(3)$ it follows that $|p|(F')=0$ if and
only if $|p^g|(F')=0$, since $m\in Q(X,{\bf C})$, consequently,
$p^g$ is equivalent to $p$ on $wX$ for each $g\in G$, where $|p|$
denotes the variation of the measure $p$. That is $p\in Q(X,{\bf
C})$. \par There exist functions $f_m$ and $f_p$ so  that
$m(dx)=f_m(x)|m|(dx)$ and $p(dy)=f_p(y)|p|(dy)$. If the measures $m$
and $p$ are $\sigma $-smooth, then their variations $|m|$ and $|p|$
are $\sigma $-smooth as well. From the construction above it follows
that the condition $|m|(X)<\infty $ implies $|p|(wX)<\infty $. In
view of Theorem I.3.2.2 \cite{bogachmtb} the Radon-Nikodym
derivatives $|m^g|(dx)/|m|(dx)$ and $|p^g|(dy)/|p|(dy)$ exist for
every $g\in G$, $x\in X$ and $y\in wX$. This implies that the
Radon-Nikodym derivatives $d_m(g,x)=m^g(dx)/ m(dx)$ and
$d_p(g,y)=p^g(dy)/p(dy)$ exist as well for all $g\in G$, $x\in X$
and $y\in wX$ (see \S I.3.2.2 \cite{bogachmtb}) such that
$$(4)\quad p^g(F')=\int_{F'} d_p(g,y)p(dy) = \lim_{ \{ F \} } m^g(F)=
 \lim_{ \{ F \} } \int_F
d_m(g,x)m(dx)$$ for each $F'$ closed in $wX$. In view of Theorem
3.6.21 \cite{eng} the topological space $wX$ is $T_1$, since $X$ is
$T_1$. On the other hand, for each point $y\in wX$ there exists an
ultrafilter of closed sets ${F'}_k$ in $wX$ satisfying the condition
$ \{ y \} = \bigcap_k {F'}_k$, since in $wX$ each singleton $ \{ y
\} $ is closed if and only if $wX$ is $T_1$ (see \S 1.5 \cite{eng}).
Hence from Formulas $(3)$ and $(4)$ of this subsection it follows
that up to $p$-almost everywhere the quasi-invariance factor
$d_p(g,y)$ can be chosen in the $y$ variable so that
$d_p(g,x)=d_m(g,x)$ for each $g\in G$ and $m$-almost everywhere on
$X$, since $X\subset wX$.

\par {\bf 12. Theorem.} {\it If $X$ is a metrizable space,
then a topological space $(Q_{\tau }^+(X), \tau _s(G))$ is
metrizable.}
\par {\bf Proof.} For each $m\in Q_{\tau }^+(X)$ put $T(m) = (m^g: g
\in G )$ and hence $T(m)\in (M_{\tau }^+(X))^G$. Since $m^e=m$,
where $e$ denotes the neutral element in the group $G$, then $T:
Q_{\tau }^+(X)\to (M_{\tau }^+(X))^G$ is the injective mapping.
\par In view of Theorem II.4.13 \cite{vardmsb61} the topological space
$(M_{\tau }^+(X), \tau _w)$ is metrizable. Let $D$ denote a metric
on $(M_{\tau }^+(X), \tau _w)$. On $(T( Q_{\tau }^+(X)))^2$ we
introduce the function
\par $(1)$ $D^G(T(m),T(p)) : = \sup_{g\in G} D(m^g,p^g)$ \\ which induces the
function \par $(2)$ $E(m,p) :=D^G(T(m),T(p)) $ \\ on $(Q_{\tau
}^+(X))^2$ for each $m, p \in Q_{\tau }^+(X)$. By the conditions of
\S 1 the mapping $h_g: X\to X$ is such that $h_g({\cal F}) \subset
{\cal F}$ for each $g\in G$, consequently, $0\le m^g(U)\le m(X)$ for
each $g\in G$ and $U\in {\cal F}$. Therefore, $E(m,p)\le m(X)+p(X)$.
On the other hand, from Formulas $(1)$ and $(2)$ of this subsection
it follows that $E(m,p)\ge D(m,p)$ for each $m, p\in Q_{\tau
}^+(X)$, consequently, $E(m,p)=0$ if and only if $m=p$. Moreover,
$E(m,p)=E(p,m)$, since $D(m^g,p^g)=D(p^g,m^g)$ for every $g\in G$
and $m, p \in Q_{\tau }^+(X)$. Furthermore, we deduce that \par
$E(m,p)=\sup_{g\in G} D(m^g,p^g)\le \sup_{g\in G} [D(m^g,s^g) +
D(s^g,p^g)]$\par $ \le E(m,s)+E(s,p)$ \\ for all $m, p, s\in Q_{\tau
}^+(X)$. Thus $E$ is the metric on $Q_{\tau }^+(X)$. \par If $(m_a:
a \in A)$ is a net in $Q_{\tau }^+(X)$ converging to $m\in Q_{\tau
}^+(X)$ relative to the $\tau _s(G)$ topology, where $A$ is a
directed set, then the net $(m_a^g: a \in A)$ converges to $m^g$
uniformly in the variable $g\in G$, hence $\lim_a m_a^g =m^g$
relative to the $\tau _w$ topology for each $g$ and uniformly in the
variable $g\in G$. This implies that $\lim_a D(m^g_a,m^g)=0$
uniformly in $g\in G$ and consequently, $\lim_a E(m_a,m)=0$.
\par Vice versa if $\lim_a E(m_a,m)=0$ for a net $(m_a: ~ a \in A)$
in $Q_{\tau }^+(X)$, then $\lim_a D(m_a^g,m^g)=0$ uniformly in the
variable $g\in G$. That is $m_a^g$ tends to $m^g$ relative to the
$\tau _w$ topology and uniformly in $g\in G$, consequently, $m_a$
tends to $m$ in $Q_{\tau }^+(X)$ relative to the $\tau _s(G)$
topology due to Theorem 2. Thus the metric $E$ on $Q_{\tau }^+(X)$
induces the equivalent topology with $\tau _s(G)$.

\par {\bf 13. Lemma.} {\it If a group $G$ acts
continuously on a topological space $X$ so that $Gx$ is dense in $X$
for each $x\in X$, $m\in Q^+(X)$, $m(A)\in \{ 0 , 1 \} $ for each
$A\in \cal F$, then $m(A)=0$ for each $A\in \cal F$.}
\par {\bf Proof.} For arbitrary two points $x\ne y\in X$ take open neighborhoods
$U_x$ and $U_y$ of $x$ and $y$ correspondingly which do not
intersect, $U_x\cap U_y=\emptyset $, since $(X\in T_1\cap
T_{3.5})\Rightarrow (X\in T_2)$ (see \cite{eng}). If $m(U_x)=1$,
then $m(U_y)=1$, since an element $g\in G$ exists so that
$(gU_x)\cap U_y\ne \emptyset $ and the measure $m$ is
quasi-invariant and non-negative. Therefore, $m(X)\ge 2$, that
contradicts the suppositions of this lemma, consequently, $m(A)=0$
for each $A\in \cal F$.

\par {\bf 14. Proposition.} {\it If a group $G$ acts
continuously on a topological space $X$ so that $Gx$ is dense in $X$
for each $x\in X$, $m\in Q^+_{\sigma }(X)$, the range $\{ m(A): A
\in {\cal F} \} =:T$ is discrete in $ [0, \infty )$ and $m(X)<\infty
$, the topological weight of $X$ is $w(X)\ge \aleph _0$, then
$m(A)=0$ for each $A\in \cal F$.}
\par {\bf Proof.} Take an arbitrary sequence $U_k$ of open subsets in $X$ so
that $U_{k+1}\subset (X\setminus (U_1\cup...\cup U_k))$ for each
$k\in \bf N$. Put $t := \inf \{ b-a: ~b>a\in T \} $. By the
conditions of this proposition $t>0$. If $m(U_k)>0$, then $m(U_l)>0$
for each $l$, since the measure $m$ is quasi-invariant and
non-negative and $Gx$ is dense in $X$ for each $x\in X$ and hence
elements $g_l\in G$ exist so that $(g_lU_l)\cap U_k\ne \emptyset $
for each $l$. Then $$m(X)\ge \sum _{k=1}^{\infty }m(U_k)=\infty ,$$
since $m(U_k)\ge t>0$ for each $k$. This contradicts the
suppositions of this proposition, hence $m(A)=0$ for each $A\in \cal
F$.

\par {\bf 15. Theorem.} {\it Suppose that $(Q^+(X)\setminus \{ 0 \} )\ne \emptyset $
and $GU=X$ for each open subset $U$ in $X$. Then the family $Q^+(X)$
is dense in the topological space $(M^+(X), \tau _w)$ relative to
the weak topology $\tau _w$.}
\par {\bf Proof.} Let $m$ be a non-trivial non-negative quasi-invariant
measure $m\in (Q^+(X)\setminus \{ 0 \} )$. If $m(V)=0$ for some
$V\in \cal F$, then $m^g(V)=0$ for each $g\in G$, since $m(V)\ge 0$
and $m$ is quasi-invariant relative to the group $G$. Therefore, in
view of Theorems I.2 and I.5 \cite{vardmsb61} (or see these theorems
in \cite{alexmsb941}) $m(U)>0$ for each open $U$ in $X$, since
$GU=X$. Since $X\in T_1\cap T_{3.5}$, for any marked point $z\in X$
a net of its neighborhoods $U_b$ exists such that $\bigcap_{b\in A}
U_b = \{ z \} $ and $U_a\subset U_b$ for each $a>b \in A$, where $A$
is a directed set. Choose a net of non-negative continuous bounded
functions $f_a$ so that $f_a(x)\le f_b(x)\le f_b(z)$ for each $a>b
\in A$ and $x\in X$, also $$\int_Xf_a(x) m(dx)=1,$$  with the
support $supp (f_b)\subset U_b$ and $f_b(z)>0$ for each $b\in A$.
Then we infer that $$\lim_a f_a(x)m(dx) = \delta _z(dx)$$ relative
to the weak topology $\tau _w$ on $M^+(X)$, where $\delta _z(dx)$
denotes the point measure in $M^+(X)$ with the support $supp (\delta
_z) = \{ z \} $ and $\delta _z(V)=1$ for each $z\in V\in \cal F$,
since $Q^+(X)\subset M^+(X)$.
\par In virtue of Theorem II.3.10 \cite{vardmsb61} the real span of point measures is dense
in $(M^+(X), \tau _w)$, consequently, $Q^+(X)$ is dense in $(M^+(X),
\tau _w)$.

\par {\bf 16. Proposition.} {\it Let a topological space
$X$ be infinite and let $X\in T_1\cap T_{3.5}$,
also $GU=X$ for each open subset $U$ in $X$, also let $Gz$ be dense
in $X$ for some $z\in X$. Then $Q^+(X)$ is not dense in $M^+(X)$
relative to the $\tau _s(G)$ topology.}
\par {\bf Proof.} We consider $(M(X),\tau _s(G))$ (see Theorem 2 above).
 Take $z\in X$ and
$\delta _z\in M^+(X)$. In the case $Q^+(X) = \{ 0 \} $ evidently
$Q^+(X)$ is not dense in $M^+(X)$ relative to the $\tau _s(G)$
topology. \par Now let $m\in (Q^+(X)\setminus \{ 0 \} )\ne \emptyset
$. From \S 15 we infer that there exists a bounded continuous
non-negative function $f$ so that $f(z)=1$, $supp (f)\subset U_b$
for some neighborhood $U_b$ of $z$ and $$\sup_{g\in G} |\int_X
f(x)(m^g(dx)- \delta _z^g(dx))|>1/3 ,$$ since $\delta ^g_z(g \{ z \}
)=1$ and an element $g\in G$ exists with $gz \in X\setminus U_b$.

\par {\bf 17. Definition.} For a positive quasi-invariant measure $m \in
Q^+(X )$ we say that its quasi-invariance factor $d_m(g,x)$ is
unbounded on $X$ for some element $g\in G$ if and only if $m ( \{ x:
x \in X; d_m(g,x)>t \} )>0$ for each $t>0$.

\par {\bf 18. Theorem.} {\it  $(1)$. The family $Q^+(X)$ induces on the
algebra $\cal F$ an uniform space structure $\cal U$ on which $G$
acts injectively. $(2)$. If there is a non-negative non-trivial
quasi-invariant measure $m\in Q^+(X)\setminus \{ 0 \}$ and a
quasi-invariance factor $d_m(g,x)$ of $m$ exists on $G\times X$ and
is unbounded on $X$ for some group element $g\in G$, then this
element $g$ acts non uniformly continuously on $\cal U.$}
\par {\bf Proof.} $(1)$. To each quasi-invariant non-negative measure
$p\in Q^+(X)$ the pseudo-metric $\rho _p(B,E) := p(B\triangle E)$ on
$\cal F$ corresponds. It induces the equivalence relation $B\Phi
_pE$ if and only if $\rho _p(B,E)=0$. Therefore, the quotient space
${\cal F}_p := {\cal F}/\Phi _p$ is the metric space. \par If a
measure $p$ is absolutely continuous relative to $m$, that is
$p<<m$, then from $B\Phi _mE$ it follows that $B\Phi _pE$. Hence
there exists the quotient mapping $\psi ^p_m : {\cal F}_p\to {\cal
F}_m$. Put $p\succeq m$ if and only if there exists a constant
$b_{p,m}>0$ so that $p(E)\le b_{p,m} m(E)$ for all $E\in {\cal F}$.
This makes on $Q^+(X)$ the directed set structure and defines the
inverse mapping system $ \{ {\cal F}_p, \psi ^p_m, Q^+(X) \} $. Put
$${\cal U} : = \lim \{ {\cal F}_p, \psi ^p_m, Q^+(X) \} ,$$
consequently, ${\cal U}$ is the uniform space. \par Since $m^g$ is
equivalent to $m$ for each $m\in Q^+(X)$ and $g\in G$, then $\rho
_{m^g}(B,E)>0$ if and only if $\rho _m(B,E)>0$, where $B, E\in \cal
F$ are arbitrary. This implies that each element $g$ of the group
$G$ acts injectively on the uniform space $\cal U$.
\par $(2)$. The condition $m\in Q^+(X)$ implies that $m^g\in Q^+(X)$ for each $g\in G$, hence
$d_m(g,x)\ge 0$ for all $x\in X$ and $g\in G$. If $m\in
Q^+(X)\setminus \{ 0 \}$ and a quasi-invariance factor $d_m(g,x)$ of
$m$ exists on $G\times X$ and is unbounded on $X$ for some $g\in G$,
then $X$ and $\cal F$ are infinite and for each $t>0$ there exist
$B, E \in {\cal F}$ so that $m^g(B\triangle E)\ge t m(B\triangle
E)$, consequently, such element $g$ acts non uniformly continuously
on $\cal U$ in this case.

\par {\bf 19. Theorem.} {\it If $X\in T_1\cap T_{3.5}$, each element $s$
in a group $G$ acts by homeomorphic and injective mappings $h_s$ on
$X$ so that $Gx$ is dense in $X$ for each $x\in X$, also $m\in
Q^+(X)\setminus \{ 0 \}$ and a quasi-invariance factor $d_m(g,x)$ of
$m$ exists on $G\times X$ and is unbounded on $X$ for some group
element $g\in G$, then this element $g$ acts discontinuously on
$\cal U.$}
\par {\bf Proof.} The imposed conditions of this proposition imply
that for each pair of points $x\ne y\in X$ and their open
neighborhoods $U_x$ and $U_y$, $x\in U_x$ and $y\in U_y$, there
exist group elements $s, q\in G$ such that $(sU_x)\cap U_y\ne
\emptyset $ and $(qU_y)\cap U_x\ne \emptyset $. In view of
Proposition 14 for every $b>0$ and $x\ne z\in X$ there exist
non-intersecting neighborhoods $U_x$ and $U_z$ so that $0<m(U_x)<b$
and $0<m(U_z)<b$. Therefore, for every $t>0$, $b>0$ and $B \in \cal
F$ there exists $E\in \cal F$ with $m(E)>0$ and $0<m(E\triangle
B)<t$ so that $m^g(E\triangle B)>bt$, consequently, $g$ induces the
mapping on $\cal U$ which is discontinuous at each element of $\cal
U$.

\par {\bf 20. Definition.} A sequence $ \{ E_k \} $ in $\cal F$ we call a
$(k',G)$ sequence if the following three conditions are satisfied:
\par $(1)$ $gE_k\uparrow X$ for each $g\in G$;
\par $(2)$ for every $k\in \bf N$, $g\in G$ and $m\in M_{\sigma
}(X)$ the sequence
\par $|m^g|_*(X\setminus E_k) := \sup \{ |m^g|(U): ~ U\subset
X\setminus E_k, ~ U \mbox{ is open} \} $ \\ converges to zero as $k$
tends to the infinity;
\par $(3)$ a function $h$ is continuous on $X$ if and only if $h$ is
continuous on $gE_k$ for every $k\in \bf N$ and $g\in G$.

\par {\bf 21. Theorem.} {\it Suppose that $(1)$ $X\in T_1\cap
T_{3.5}$; \par $(2)$ each element $s$ in the group $G$ acts by
homeomorphic and injective mappings $h_s$ on $X$; \par $(3)$ a
sequence $\{ m_k : k \in {\bf N} \} \subset Q_{\sigma }(X)$
converges to $m_0\in M(X)$ relative to the $\tau _w$ topology, also
\par $(4)$ a function $p: G\times X \to [0,\infty )$ exists such that
for each $n\in {\bf N}$ and $g\in G$ the inequality $ ~
|d_{m_n}(g,x)|\le p(g,x)$ is fulfilled $m_n$-a.e. in the $x$
variable, where $$p(g,*)\in \bigcap_{n=0}^{\infty } L^1(m_n)\mbox{
and  } \sup_{0\le n} \| p(g,*) \| _{L^1(m_n)}< \infty .\mbox{
Then}$$
\par $(5)$ for each $(k',G)$ sequence $ \{ E_k: k \in {\bf N} \} $
and each element $g\in G$ the limit \par $~ \lim_{k\to \infty } |
m_n^g|_*(X-E_k) =0 $ converges uniformly in $n$ and \par $(6)$ $m_0
\in Q_{\sigma }(X)$.}
\par {\bf Proof.} Consider a $(k',G)$ sequence $ \{ E_k: k \in {\bf N}
\} $. Put $$(7)\quad D_g(h,s) := \sum_{k=1}^{\infty } 2^{-k}
\sup_{x\in gE_k} |h(x)-s(x)|$$ for each $h, s \in S_1$, where $$ S_1
:= \{ h: ~ h \in C_b(X,{\bf R}), \| h \| \le 1 \} ,$$ $g$ is a
marked element in $G$. Then from Conditions $(1)$ and $(2)$ and
Formula $(7)$ of this subsection and Conditions 20$(1-3)$ of
Definition 20 we infer that the mapping $D_g$ is the metric on $S_1$
and the metric space $(S_1,D_g)$ is complete.
\par Choose an arbitrary quasi-invariant $\sigma $-smooth measure
$m\in Q_{\sigma }(X)$ and define the functional $$F(g,h) :=
\int_Xh(x)m^g(dx)$$ on $S_1$, where $h\in S_1$. Applying
Radon-Nikodym's theorem to equivalent measures  $m^g \sim m$ in
$M_{\sigma }(X)$ one gets the densities $d_m(g,x)=m^g(dx)/m(dx)$.
From the identity $$\int_Xh(x)m^g(dx)=\int_Xh(gy)m(dy)$$ and
Conditions 20$(1-3)$ of Definition 20 we deduce that the functional
$F(g,*)$ is continuous on the metric space $(S_1,D_g)$. Therefore,
for every $b>0$ and $m\in Q_{\sigma }(X)$ \par the set $ \{ h: ~
h\in S_1, ~ |\int_X h(x)m^g(dx) | \le b \} $ is closed in
$(S_1,D_g)$. \par Then each set of the form
$$W_k(b)=W_k^g(b) := \{ h: ~ h \in S_1, \sup_{l, n \ge k}
|\int_X h(x)m_l^g(dx) - \int_Xh(x)m_n^g(dx)|\le b \} $$ is closed in
$(S_1,D_g)$ and $$S_1=\bigcup_{k=1}^{\infty }W_k(b)\mbox{ for each }
0<b<1.$$ In view of the Baire category theorem 3.9.3 \cite{eng} a
natural number $u$ exists so that $W_u(b)$ contains an open subset
in $(S_1,D_g)$. \par Therefore, there exist $u\in \bf N$, $h_0\in
S_1$, $j\in \bf N$ and $t>0$ such that from $n, k \ge u$, $h\in S_1$
and $|h(x)-h_0(x)|<t$ on $E_j$ it follows that $$ |\int_X
h(x)m_n^g(dx) - \int_X h(x)m_k^g(x)|\le t,\mbox{  particularly,}$$
$$ |\int_{X-E_j}m_n^g(dx)-\int_{X-E_j}m_k^g(dx)|\le 2t $$ if take
$h_0$ satisfying the restriction $h_0|_{E_j}=0$. In virtue of Lemma
I.3 \cite{vardmsb61} and Conditions 20$(1-3)$ of Definition 20 the
estimate $$|m_n^g-m_k^g|_*(X-E_j)\le 2t$$ is fulfilled for all $n,k
\ge u$. Choose a natural number $i(1)\ge j$ so that
$$|m_u^g|(X-E_{i(1)})\le t,\mbox{ consequently, }$$ $$|m_n^g|_*(X-E_{i(1)})\le
[|m_n^g-m_j^g|_*(X-E_{i(1)})+|m_j^g|_*(X-E_{i(1)})]\le 3t.$$ Then
one can take a natural number $i(2)\ge i(1)$ for which $$\sup_{l<u}
|m_l^g|_*(X-E_{i(2)})\le 3t,$$ consequently, for each $b>0$ there
exists $i(2)$ so that $$\sup_{ n\in \bf N}|m_n^g|_*(X-E_{i(2)}) \le
3b.$$ This implies assertion $(5)$.
\par From the convergence of the sequence $m_n$ to $m_0$ relative to the $\tau
_w$ topology, Conditions $(1)$ and $(2)$ of this subsection and
Conditions 20$(1-3)$ of Definition 20 it follows that
$$\underline{\lim}_{n\to \infty } |m_n^g|_*(X-E_i)\ge
|m_0^g|_*(X-E_i)$$ for each $g\in G$, consequently,
$$|m_0^g|_*(X-E_i)\le 3b$$ for all $i\ge i(2)$ and an arbitrary marked
$g$ in $G$ and $E_i$, $b$, $i(2)$ described above in this
subsection. That is for each $g\in G$ the limit $$\lim_{i\to \infty
} |m_i^g|_*(X-E_i)=0$$ exists for each $(k',G)$ sequence $ \{ E_i: i
\in {\bf N} \} $.
\par We consider a sequence $ \{ Z_n: n \in {\bf N} \} $ of closed
in $X$ subsets for which the following conditions are valid:
\par $(8)$ for each $g\in G$ $~gZ_n\uparrow X$ as $n\to \infty $;
\par $(9)$ for each $n\in \bf N$ there exists an open subset $U_n$ so
that $gZ_n\subset gU_n\subset gZ_{n+1}$ for each $g\in G$. \par From
Theorem I.12 \cite{vardmsb61} and Conditions $(1)$ and $(2)$ it
follows that if a sequence $ \{ Z_n: n \in {\bf N} \} $ of closed in
$X$ subsets satisfies Conditions $(8)$ and $(9)$, then it is a
$(k',G)$ sequence, since $h_s: X\to X$ is the homeomorphism for each
$s\in G$, $h_s(x)=sx$. Hence $$\lim_{i\to \infty } |m_0^g|(X-Z_i)
=0$$ for each $g\in G$. Then by Theorem I.19 \cite{vardmsb61}
$m_0^g\in M_{\sigma }(X)$ for each $g\in G$. \par Let now $V\in \cal
F$ and $|m_0|(V)=0$. Then we get that $$\lim_{n\to \infty } |m_n|(V)
=0.$$ There are valid the identity
$$(10)\quad \int_Xf(x)m_n^g(dx) = \int_X
f(x)d_{m_n}(g,x)m_n(dx) \quad \mbox{ and the inequality} $$
$$(11)\quad |\int_Xf(x)d_{m_n}(g,x)m_n(dx)|\le \int_X p(g,x)
|f(x)||m_n|(dx)$$ for each $g\in G$ and $n\in \bf N$. For each open
subset $U$ in the topological space $X$, $g\in G$ and $k\in \bf N$
the inequality $$|m_k^g|(U) \le
 |m_k|(U) \sup_{0\le n} \| p(g,*) \|_{L^1(m_n)}$$ is fulfilled.
  Take a sequence $f_u$ of continuous functions $f_u: X\to [0,1]$
so that $f_u(x)=1$ for each $x\in V$ and $u\in \bf N$ and $$\int_X
f_u(x)|m_0|(dx)<b_u,$$ where $b_u\downarrow 0$, consequently,
$$\lim_u \int_Xf_u(x)|m_0|(dx)=|m_0|(V)=0$$
(see \cite{alexmsb941,bogachmtb,vardmsb61}). Applying Condition
$(4)$ we infer from Formulas $(10)$ and $(11)$ of this subsection
that for each $b>0$ and $g\in G$ there exist natural numbers $n$ and
$v$ such that $$\int_X f_v(x)|m_l^g|(dx) <b\mbox{ for each }l>n.$$
On the other hand, the inequality is valid: $$|m_l^g|(V)\le \int_X
f_v(x)|m_l^g|(dx).$$ Therefore, the limit $$\lim_{n\to \infty }
|m_n^g|(V)=0$$ exists for each $g\in G$, consequently,
$|m_0^g|(V)=0$ for each $g\in G$. That is the measure $m^g$ is
equivalent with the measure $m_0$ for each group element $g\in G$.
Thus the measure $m_0$ is quasi-invariant and consequently, $m_0\in
Q_{\sigma }(X)$, since $Q(X)\cap M_{\sigma }(X)=Q_{\sigma }(X)$.

\par {\bf 22. Theorem.} {\it Let $G$ be a subgroup of a topological group
$G_1$ and let Conditions 21$(1)$ and $(2)$ of subsection 21 be
satisfied for $G$ and $G_1$, let also the topological character of
$G_1$ be $ \chi (G_1)=\aleph _0$ and let the mapping $G\times X\ni
(g,x)\mapsto gx\in X$ be continuous. Suppose also that for a
quasi-invariant measure $m\in Q_{\sigma }(X,{\bf F},{\cal F},G)$ a
sequence $m_k:=m^{g_k}$ fulfills Conditions 21$(3)$ and $(4)$ of
subsection 21 relative to $G$ for each $g_0\in G_1$ and each
sequence $g_k$ in $G$ converging to $g_0$ relative to the left
uniformity induced by the topology on $G_1$. Then $m\in Q_{\sigma
}(X,{\bf F},{\cal F},G_1)$.}
\par {\bf Proof.} For the topological group $G_1$ there exists the
left uniformity on $G_1$ induced by the topology on $G_1$ (see \S
8.1.17 in \cite{eng}). Since $ \chi (G_1)=\aleph _0$, then for each
$g_0\in G_1$ there is a sequence $g_k$ in $G$ converging to $g_0$
relative to the left uniformity of $G_1$. In view of Theorem 21 the
limit exists $$\lim_k m_k=m_0\in Q_{\sigma }(X,{\bf F},{\cal
F},G),$$ we denote it by $m_0 =: m^{g_0}$ also, since $$\lim_k
g_k=g_0\mbox{ and } m_k:=m^{g_k}.$$
\par Suppose that there are two sequences $s_k$ and $q_k$ converging
to $g_0$. Put $g_{2k}=s_k$ and $g_{2k+1}=q_k$ for each $k$, hence
$g_k$ converges to $g_0$ as well and consequently, $$\lim_k
m^{s_k}=\lim_k m^{q_k}=\lim_k m^{g_k},$$ since the limit $\lim_k
m^{g_k}$ exists according to Condition 21$(3)$ of subsection 21.
Thus $m^{g_0}$ is independent of the sequence $g_k$ in $G$
converging to $g_0\in G_1$.
\par  Since the mapping $G\times X\ni (g,x)\mapsto gx\in X$ is
continuous and $h_g$ acts by homeomorphisms on $X$ for each $g\in
G_1$ and $$\lim_k g_k=g_0,\mbox{ then }g_0U = \bigcup_{n=1}^{\infty
} \bigcap_{k=n}^{\infty } g_kU$$ for each $U$ open in $X$.
Therefore, $$|m|(g_0U)= \lim_n |m|(\bigcap_{k=n}^{\infty } g_kU)$$
for each $m\in M_{\sigma }(X)$. If $|m|(V)=0$ for some $V\in \cal F$
and $m\in M_{\sigma }(X)$, then for each $b>0$ there exists an open
subset $K_b$ in $X$ such that \par $V\subset K_b$ and $|m|(K_b)<b$
\\  \cite{alexmsb941,bogachmtb,vardmsb61}. If $m\in Q_{\sigma }(X,{\bf F},{\cal F},G)$,
then $|m^g|(V)=0$ for each $g\in G$ as well. On the other hand,
$$|m|(\bigcap_{k=n}^{\infty } g_kU)\le |m| (g_kU)$$ for each $k\ge n$.
Then from Condition 21$(4)$ and Formulas 21$(10)$ and $(11)$ of
subsection 21 we infer that $|m^{g_0}|(V)=0$. That is these measures
are equivalent $m \sim m^{g_0}$ for each $g_0\in G_1$. Thus $m\in
Q_{\sigma }(X,{\bf F},{\cal F},G_1)$.

\par {\bf 23. Theorem.} {\it Let $ \{ G_i, h^i_j, J \} $ and $ \{ X_i, p^i_j, J \} $
be inverse mapping systems of groups $G_j$ and topological spaces
$X_j$ so that $G_j: X_j\to X_j$, where $h^i_j: G_i\to G_j$ are
homomorphisms of groups and $p^i_j: X_i\to X_j$ are continuous
mappings for each $i\ge j$ in a directed set $J$. Let also a
topological space $X$ be homeomorphic with $\lim \{ X_j, p^i_j, J \}
$. Then each measure $m\in Q(X,{\bf F},{\cal F},G)$ has the
decomposition $$m=\lim \{ m_i; p^i_j; J \} $$ with $m_i\in
Q(X_i,{\bf F},{\cal F}_i,G_i)$ for each $i\in J$, which induces the
continuous mapping from $Q(X,{\bf F},{\cal F},G)$ into $ \lim \{
Q(X_i,{\bf F},{\cal F}_i,G_i); p^i_j; h^i_j; J \} $ relative to
their $\tau _w(G)$ and $\{ \tau _w(G_i); p^i_j; h^i_j; J \} $
topologies. }
\par {\bf Proof.} Consider continuous projective mappings $p^i: X\to X_i$
and projective group homomorphisms $h^i: G \to G_i$ for each $i\in
J$ corresponding to these inverse mapping systems. \par If $m\in
Q(X,{\bf F},{\cal F},G)$, then $$m((p^i)^{-1}(A)) =: m_i(A)$$ is a
measure on ${\cal F}_i$, since $(p^i)^{-1}({\cal F}_i)\subset \cal
F$ for each $i\in J$. On the other hand, $g((p^i)^{-1}(A))
=(p^i)^{-1}(g_iA)$ for each $A\in {\cal F}_i$ and $g\in G$, where
$g_i=h^i(g)$, since $G_i: X_i\to X_i$, $~p^i_j\circ p_i=p_j$ and
$h^i_j\circ h_i=h_j$ for each $i\ge j\in J$. This family of measures
forms and inverse system, since $$m_i((p^i_j)^{-1}(B)) =
m((p^j)^{-1}(B)) = m_j(B)$$ for each $i>j \in J$ and $B\in {\cal
F}_j$. If $|m|_i(A) =0$ for some $A\in {\cal F}_i$, then
$|m|((p^i)^{-1}(A)) =0$, consequently, $|m^g|((p^i)^{-1}(A)) =0$ for
each $g\in G$ and hence $|m_i^{g_i}|(A)=0$. Thus the measure $m_i$
on ${\cal F}_i$ is quasi-invariant relative to the group $G_i$ for
each $i\in J$. Then there is valid the equality $$(1)\quad \int_X
S(f_i)(y)m(dy) = \int_{X_i} f_i(x_i)m_i(dx_i)$$ for each $f_i\in
C_b(X_i,{\bf R})$, where $S(f_i)$ is the cylinder function
corresponding to $f_i$, $S(f_i)(y)=f_i(p_i(y))$ for each $y\in X$.
In view of Proposition 2.5.5 \cite{eng}, Theorem 2 and Formula $(1)$
of this subsection this mapping $ m\mapsto \lim \{ m_i; p^i_j; J \}
$ is continuous from $Q(X,{\bf F},{\cal F},G)$ into $ \lim \{
Q(X_i,{\bf F},{\cal F}_i,G_i); p^i_j; h^i_j; J \} $ relative to
their $\tau _w(G)$ and $\{ \tau _w(G_i); p^i_j; h^i_j; J \} $
topologies.

\par {\bf 24. Theorem.} {\it Let suppositions of Theorem 23 be
satisfied, where $PQ^+_{\sigma }(X_i,{\bf F},{\cal B}_i,G_i) := \{
m_i: ~ m_i\in Q^+_{\sigma }(X_i,{\bf F},{\cal B}_i,G_i), m_i(X_i)=1
\} $ is on the algebra ${\cal B}_i = {\cal B}_i(X_i)$ for each $i\in
J$. Let also $(X_i,{\cal B}_i)$ be a Radon space for each $i\in J$.
If $J$ is countable, then there exists an embedding of $PQ^+_{\sigma
}(X,{\bf F},{\cal B},G)$ into $\lim \{ PQ^+_{\sigma }(X_i,{\bf
F},{\cal B}_i, G_i); h^i_j; p^i_j; J \} $ relative to the weak $\tau
_w(G)$ and  $\{ \tau _w(G_i); p^i_j; h^i_j; J \} $ topologies.}
\par {\bf Proof.} Each $\sigma $-smooth measure $m$ defined on $\cal F$
has the natural extension on $\cal B$ and analogously for each $m_i$
\cite{alexmsb941,bogachmtb,vardmsb61}. In view of Theorem 23 if
$m\in PQ^+_{\sigma }(X,{\bf F},{\cal B},G)$, then $ \{ m_i; p^i_j; J
\} \in \{ PQ^+_{\sigma }(X_i,{\bf F},{\cal B}_i,G_i); h^i_j; p^i_j;
J \} $, since the conditions \par $A_{i,k}\in {\cal B}_i$ and
$A_{i,k}\cap A_{i,n}=\emptyset $ for each $k\ne n \in \bf N$ imply
that \par $(p^i)^{-1}(A_{i,k})\in \cal B$ and
$(p^i)^{-1}(A_{i,k})\cap (p^i)^{-1}(A_{i,n})=\emptyset $ for each
$k\ne n \in \bf N$. \par Suppose that the inverse system of measures
$\{ m_i; p^i_j; J \} $ is given so that $m_i\in PQ^+_{\sigma
}(X_i,{\bf F},{\cal B}_i, G_i)$ for each $i\in J$. This system
induces the quasi-measure $q$ on $${\cal H} := \bigcup_{i\in J}
(p^i)^{-1}({\cal B}_i)\mbox{ so that}$$
 $m_i((p^i_j)^{-1}(B)) = q((p^j)^{-1}(B)) = m_j(B)$ for each
$i>j \in J$ and $B\in {\cal B}_j$. \par Moreover, $q(X)=1$, since
$m_i(X_i)=1$ for each $i\in J$. The algebra of cylinder subsets
${\cal H}$ is contained in the algebra ${\cal B}$, since
$(p^i)^{-1}({\sf Z}_i)\subset \sf Z$ and $(p^i)^{-1}({\sf
U}_i)\subset \sf U$ for each $i\in J$. Thus $q$ is a bounded
non-negative quasi-measure. In view of Kolmogorov's Theorem I.1.3
and Propositions I.1.7 and I.1.8 \cite{dal} (or see
\cite{vahtarchobb}) it has the unique $\sigma $-smooth extension $q$
on the minimal $\sigma $-algebra $\sigma {\cal H}$ containing ${\cal
H}$. \par The base of the topology of $X$ consists of open subsets
of the form $\bigcap_{l=1}^n(p^{i(l)})^{-1}(A_{i(l)})$ with
$A_{i(l)}\in {\sf U}_{i(l)}$ for each $i(l)\in J$, $l=1,...,n$ and
$n\in \bf N$. Since $J$ is countable, one gets that ${\cal B}=\sigma
{\cal H}$. Thus $q$ and $m$ coincide on $\cal B$, if $\{ m_i; p^i_j;
J \} $ is constructed from $m\in PQ^+_{\sigma }(X,{\bf F},{\cal
B},G)$. That is the mapping from $PQ^+_{\sigma }(X,{\bf F},{\cal
B},G)$ into $\lim \{ PQ^+_{\sigma }(X_i,{\bf F},{\cal B}_i, G_i);
h^i_j; p^i_j; J \} $ is injective.
\par Therefore, for each $f\in C_b(X,[0,\infty ))$ there exists a
net $f_i\in C_b(X_i,[0,\infty ))$ so that
$$(1)\quad \int_X f(x)m(dx) = \lim_i \int_{X_i}
f_i(x_i)m_i(dx_i),\mbox{  where }$$ $$(2)\quad f((p^i)^{-1}(x_i))\le
f_i(x_i)$$ for each $x_i\in X_i$ and $i\in J$. Particularly, if $f:
X\to [0,1]$, then one can choose \par $(3)$ $f_i: X_i\to [0,2]$ for
each $i\in J$. \par Then we deduce that
$$(4)\quad \int_Xf(x)m^g(dx)
=\lim_i \int_{X_i}f_i(x_i)m_i^{g_i}(dx_i)$$
$$=\lim_i\int_{X_i}f_i(g_iy_i)m_i(dy_i) =\int_Xf(gy)m(dy)$$ for each
$g\in G$. If $m(A)=0$ for some $A\in {\cal B}$, then $m^g(A)=0$ for
each $g\in G$. Then for each $t>0$ and $g_1\in G,...,g_n\in G$,
$n\in \bf N$, there exists $f\in C_b(X,[0,1])$ so that
$$f|_A=1\mbox{  and  }\int_Xf(x)m^{g_l}(dx)<t\mbox{  for each  }l=1,...,n,$$
consequently, there exists a net $\{ f_i: i \} $ satisfying
Conditions $(1-3)$ above and $k\in J$ so that
$$\int_Xf(x)m^{g_l}(dx)\le
\int_{X_i}f_i(x_i)m_i^{h^i(g_l)}(dx_i)<2t$$
for each $i\ge k\in J$ and $l=1,...,n$.

\par From Proposition 2.5.5 and Theorem 2.5.8 \cite{eng}, Formulas $(1-4)$
of this subsection and Theorem 2 we infer that this embedding of
$PQ^+_{\sigma }(X,{\bf F},{\cal B},G)$ into $\lim \{ PQ^+_{\sigma
}(X_i,{\bf F},{\cal B}_i, G_i); h^i_j; p^i_j; J \} $ is continuous
and open relative to the weak $\tau _w(G)$ and $\{ \tau _w(G_i);
p^i_j; h^i_j; J \} $ topologies, where $PQ^+_{\sigma }(X,{\bf
F},{\cal B}, G)$ and all $PQ^+_{\sigma }(X_i,{\bf F},{\cal B}_i,
G_i)$ are considered relative to their $\tau _w(G)$ and $\tau
_w(G_i)$ topologies.


\begin{thebibliography}{99}

\bibitem{alexmsb941} A.D. Alexandroff. {\it Additive set-functions in
abstract spaces (1-4)}// Sb. Mathem. (1) {\bf 8: 2 } (1940),
307-312; (2) {\bf 9: 3} (1941), 563-628; (3) {\bf 13: 3} (1943),
169-238.

\bibitem{beldal} Ya.I. Belopolskaya, Yu.L. Dalecky.
{\it Stochastic equations and differential geometry} (Dordrecht:
Kluwer Acad. Publ., 1989).

\bibitem{bogachmtb} V.I. Bogachev. {\it Measure theory. V. 1, 2}
(Berlin: Springer-Verlag, 2007).

\bibitem{bourm70} N. Bourbaki. {\it Int\'egration.} Livre VI.
Fasc. XIII, XXI, XXIX, XXXV. Ch. 1-9 (Paris: Hermann, 1963, 1965, 1967,
1969).

\bibitem{constanb} C. Constantinescu. {\it Spaces of measures}
(Berlin: Walter de Gruyter, 1984).

\bibitem{dal} Yu. L. Dalecky, S.V. Fomin. {\it Measures and Differential
Equations in Infinite-Dimensional Spaces} (Dordrecht: Kluwer Acad.
Publ., 1991).

\bibitem{dalshn} Yu.L. Dalecky, Ya.L. Shnaiderman.
{\it Diffusion and quasi-invariant measures on infinite-dimensional
Lie groups} // Funct. Anal. and its Applic. {\bf 3} (1969), 156-158.

\bibitem{eng} R. Engelking. {\it General topology} (Moscow: Mir, 1986).

\bibitem{ferpelruizb} J.C. Ferrando, M. L\'opez Pellicer, L.M.
S\'anches Ruiz. {\it Metrizable barrelled spaces}. Pitman Research
Notes in Mathematics Series. V. {\bf 332} (New York: John Wiley and
Sons Inc., Longman, 1995).

\bibitem{fidal} F. Fidaleo. {\it Continuity of Borel actions
of Polish groups on standard measure algebras} // Atti Sem. Mat.
Fiz. Univ. Modena {\bf 48} (2000), 79-89.

\bibitem{hew} E. Hewitt, K.A. Ross. {\it Abstract harmonic analysis.
V. 1, 2} (Berlin: Springer-Verlag, 1994).

\bibitem{kuob} H.-H. Kuo. {\it Gaussian measures in Banach spaces}
(Berlin: Springer-Verlag, 1975).

\bibitem{ludjms2005} S.V. Ludkovsky. {\it Quasi-invariant and
pseudo-differentiable measures with values in non-archimedean fields
on a non-archimedean Banach space} // J. of Mathem. Sci., N.Y.
(Springer) {\bf 128: 6} (2005), 3428-3460.

\bibitem{ludanmat2002} S.V. Ludkovsky. {\it Quasi-invariant and
pseudo-differentiable real valued measures on non-archimedean Banach
spaces}//  Analysis Mathematica {\bf 28} (2002), 287-316.

\bibitem{ludqupdmb} S.V. Ludkovsky. {\it Quasi-invariant and
pseudo-differentiable measures in Banach spaces}. ISBN
978-1-60692-734-2 (New York: Nova Science Publishers, Inc., 2009).

\bibitem{bulopnlcga:12} S.V. Ludkovsky. {\it Operators on a non
locally compact group algebra}// Bull. Sciences Math\'ematiques
(Paris). Ser.2, {\bf 137} (2013), 557-573; DOI:
10.1016/j.bulsci.2012.11.008

\bibitem{lugmj2014} S.V. Ludkovsky. {\it Meta-centralizers of non locally
compact group algebras} // Glasgow Mathematical Journal, {\bf 57} (2015), 349-364; DOI:
10.1017/S0017089514000330.

\bibitem{lujms147:3:08} S.V. Ludkovsky.
{\it Topological transformation groups of manifolds over
non-Archimedean fields, representations and quasi-invariant
measures, I, II} // J. of Mathem. Sci., N.Y. (Springer) (I) {\bf
147: 3} (2008), 6703-6846; (II) {\bf 150: 4} (2008), 2123-2223.

\bibitem{lusmldg05} S.V. Ludkovsky. {\it Stochastic processes on
geometric loop groups, diffeomorphism groups of connected manifolds,
associated unitary representations}// J. of Mathem. Sci., N.Y.
(Springer) {\bf 141: 3} (2007), 1331-1384.

\bibitem{lurimut98} S.V. Ludkovsky. {\it Irreducible unitary representations of a
diffeomorphisms group of an infinite-dimensional real manifold}
// Rendic. dell'Istituto di Matematica dell'Univ. di
Trieste. Nuova Serie, {\bf 30} (1998), 21-43.

\bibitem{lurimut99} S.V. Ludkovsky. {\it  Quasi-invariant measures on a group of
diffeomorphisms of an infinite-dimensional real manifold and induced
irreducible unitary representations} // Rendic. dell'Istituto di
Matematica dell'Univ. di Trieste. Nuova Serie, {\bf 31} (1999),
101-134.

\bibitem{ludan2000} S.V. Ludkovsky. {\it  Quasi-invariant measures on loop groups of
Riemann manifolds} //  Dokl. Mathemat. {\bf 61: 1} (2000), 54-56.

\bibitem{ludspb} S.V.Ludkovsky. {\it  Stochastic processes in non-archimedean Banach
spaces, manifolds and topological groups} ISBN   978-1-61668-787-8
(New York: Nova Science Publishers, Inc., 2010).

\bibitem{luambp99} S.V. Ludkovsky. {\it Properties of quasi-invariant
measures on topological groups and associated algebras} // Annales
Math. B. Pascal. {\bf 6: 1} (1999), 33-45.

\bibitem{skorohodb} A.V. Skorohod. {\it Integration in Hilbert
space} (Berlin: Springer-Verlag, 1974).

\bibitem{naribeckb} L. Narici, E. Beckenstein. {\it Topological vector
spaces} (Marcel Dekker, Inc.: New York, 1985).

\bibitem{smolfom} O. G. Smolyanov, S. V. Fomin. {\it Measures on linear topological
spaces}// Russ. Math. Surv, {\bf 31: 4} (1976), 3--56.

\bibitem{vahtarchobb} N. N. Vahaniya, V. I. Tarieladze,
S. A. Chobanyan. {\it Probability distributions in Banach spaces}
(Moscow: Nauka, 1985).

\bibitem{vardmsb61} V.S. Varadarajan. {\it Measures on topological spaces} // Mat. Sbornik.
{\bf 55} (1961), 35-100 (in Russian); English transl.: Amer. Math.
Soc. Transl. (2) {\bf 48} (1965), 161-228.

\end{thebibliography}
\end{document}